\documentclass{article}[10pt]
\usepackage{amsmath}
\usepackage{amsthm}
\usepackage{amssymb}
\usepackage{llncsdoc}

\newtheorem{theorem}{Theorem}[section]
\newtheorem{prop}[theorem]{Proposition}
\newtheorem{corollary}[theorem]{Corollary}
\newtheorem{definition}[theorem]{Definition}

\newtheorem{example}[theorem]{Example}
\newtheorem{lemma}[theorem]{Lemma}

\newtheorem{problem}{Problem}

\title{Non-characterizability in belief revision: \\ an application of finite model theory}
\author{Gy\"orgy Tur\'an\thanks{Partially supported by NSF grant CCF-0916708}  \\ University of Illinois at Chicago\\ MTA-SZTE Research Group on Artificial Intelligence \and Jon Yaggie$^*$\\University of Illinois at Chicago}
\begin{document}
\maketitle

\begin{abstract}
A formal framework is given for the characterizability of a class of belief revision operators, defined using minimization over a class
of partial preorders, by postulates. It is shown that for partial orders characterizability implies a definability property
of the class of partial orders in monadic second-order logic.
Based on a non-definability result for a class of partial orders, an example is given of a non-characterizable class of revision operators.
This appears to be the first non-characterizability result in belief revision.
\end{abstract}

\section{Introduction}

The main approach to belief change is the AGM approach pioneered by \cite{AGM}. It provides many characterization
results for belief change operators in terms of rationality postulates~\cite{H1999}. \emph{Are there cases where no characterization can be given?} Answering
this non-axiomatizability question presupposes a formal definition of a postulate. However, as noted in the survey paper \cite{AGM25}
\begin{quote}
``theories of belief change developed in the AGM tradition are not logics in a strict sense, but rather informal axiomatic theories of
belief change. Instead of characterizing the models of belief and belief change in a formalized object language, the AGM approach uses a natural
language (ordinary mathematical English) to characterize the mathematical structures under study.''
\end{quote}
Ferm\'e and Hansson then proceed to describe modal and dynamic logic approaches to modeling belief change (see, e.g., \cite{DHWKB}).
As far as we know, the question of non-characterizability by
postulates has not been considered before in those frameworks either.

In this note we provide a formal framework for studying characterizability, based on the approach of Katsuno and Mendelzon~\cite{KM}.
A revision operator $*$ in \cite{KM} is considered to assign a revised knowledge base $K * \varphi$ to every knowledge base $K$ and every revising formula $\varphi$.
However, the results remain valid if one considers $*$ to act on a fixed knowledge base $K$ and an arbitrary revising formula $\varphi$.
We are not discussing iterated revision here, and so there is no interaction between the revisions of different knowledge bases.
Katsuno and Mendelzon prove the following results.

\begin{theorem} \cite{KM} \label{th:km}

\emph{a)} There is a finite set of postulates such that a revision operator satisfies these postulates iff there is a faithful \emph{total preorder}
representing it with minimization.

\emph{b)}There is a finite set of postulates such that a revision operator satisfies these postulates iff there is a faithful \emph{partial preorder}
representing it with minimization.

\emph{c)}There is a finite set of postulates such that a revision operator satisfies these postulates iff there is a faithful \emph{partial order}
representing it with minimization.
\end{theorem}

Part \emph{a)} is a finite version of Grove's characterization of the AGM postulates in terms of systems of spheres \cite{G}.
The postulates for parts \emph{b)} and \emph{c)} are the same, and different from \emph{a)}.

We consider the following general question.

\begin{problem} \label{pr:gen}
Let ${\cal R}$ be a family of partial preorders. Is there a finite set of postulates such that a revision operator satisfies these postulates
iff there is a faithful \emph{partial preorder from ${\cal R}$} representing it with minimization?
\end{problem}


Our goal is to prove a negative answer for a particular family ${\cal R}$. The formal definition of characterizability aims at providing a formal framework
for proving this negative result.
Formulating frameworks for other types of non-characterizability results appears to be an interesting topic for future work.

Non-characterizability is proved using
a translation from postulates to universal monadic second-order formulas over the language of partial preorders.
It follows from this translation that for partial orders postulate characterization
of the class of revision operators implies universal monadic second-order definability of the class of partial orders considered.
Thus, non-definability of the class of partial orders implies non-characterizability  by postulates.
Non-definability in monadic second-order logic is a well-studied topic in finite model theory \cite{EF,L}.
We give such a non-definability result for a particular class of partial orders.
This class is constructed to give a first example of non-characterizability.
It remains an interesting problem to find more natural examples.
A candidate is the class of 2-dimensional partial orders.

\section{Preliminaries}

We consider propositional logic knowledge bases $K$ over a fixed finite set of variables.
We write $K_n$ to indicate that $K$ is over $n$ variables.
Truth assignments (or interpretations) are assignments of truth values to the variables.
The set of truth assignments satisfying a formula $\varphi$ is denoted by $|\varphi|$.
Given a set $A$ of truth assignments, $\langle A \rangle$ is some formula $\varphi$ such that $|\varphi| = A$.
A knowledge base is represented by a single formula~\footnote{We note that because of finiteness this representation corresponds to the
belief set framework. Computational complexity issues are not discussed here thus the details of the representation are irrelevant.}.

Given a knowledge base $K$, a belief revision operator $*$ assigns a formula $K * \varphi$ to every formula $\varphi$. Here $\varphi$ is called
the revising formula, and $K * \varphi$ is called the revised knowledge base.

A partial preorder is $R = (U, \le)$, where $U$ is a finite ground set and $\le$ is a reflexive, transitive binary relation.
A partial order is, in addition, antisymmetric. We write $a \sim b$ if $a$ and $b$ are incomparable.
The comparability graph of $R$ is the undirected graph over $U$ such that
for any pair of vertices $(a, b)$ is an edge iff $a \le b$ or $b \le a$.
An element $a$ is minimal if there is no $b$ such that $b < a$, where $b < a$ iff $b \le a$ but $a \not\le b$.
If $U' \subseteq U$ then $a$ is minimal in $U'$ if $a \in U'$ and there is no $b \in U'$ such that $b < a$.
The set of minimal elements of $U'$ is denoted by $\min\nolimits_\le U'$.

\begin{definition} (Faithful partial preorder)
A \emph{faithful partial preorder} for a knowledge base $K_n$ is a pair $F = (R, t)$, where $R = (U, \le)$ is a partial preorder on $2^n$ elements
and $t : U \to \{0, 1\}^n$ is a bijection between the elements of $U$ and truth assignments, such that
\begin{enumerate}
\item $a \in U$ is minimal iff $t(a)$ satisfies $K_n$,

\item if $t(a)$ satisfies $K_n$ and $t(b)$ falsifies $K_n$ then $a < b$.
\end{enumerate}
\end{definition}

In the standard definition the partial preorder is defined over the set of truth assignments. For our discussion
it is more convenient to separate the partial preorder and the labeling of its elements by truth assignments.

The basic construction used in Theorem~\ref{th:km} is that of a revision operator determined by a faithful partial preorder
using minimization.

\begin{definition} (Revision using minimization) \label{def:min}
The \emph{revision operator $*_F$ for $K$, determined by a faithful partial preorder $F$ for $K$, using minimization} is
\begin{equation} \label{eq:min} \nonumber
 K *_F \varphi = \langle \min\nolimits_{\le} t^{-1}(|\varphi|)\rangle.
\end{equation}
\end{definition}

Thus the revised knowledge base is satisfied by the minimal satisfying truth assignments of the revising formula.
Faithfulness implies that if the revising formula is consistent with the knowledge base then the revised knowledge base
is the conjunction of the knowledge base and the revising formula.

We will use some notions from finite model theory.
General introductions to the topic are given in \cite{EF, L}.
The notions used are introduced in the later sections, so our discussion is essentially self-contained.

\section{Postulates}

Consider the AGM postulates
\begin{equation} \label{eq:agm1}
\textrm{if $K$ is satisfiable then $K * \varphi$ is also satisfiable}
\end{equation}
and
\begin{equation} \label{eq:agm2}
 \textrm{if $(K * \varphi) \wedge \psi$ is satisfiable then
 $K * (\varphi \wedge \psi) \vdash (K * \varphi) \wedge \psi$.}
 \end{equation}
Here $K, K * \varphi, \varphi$ and $\psi$ can be considered as unary predicates over the set of interpretations,
and thus (\ref{eq:agm1}) can be rewritten as
\begin{equation} \label{eq:agm3} \nonumber
[\exists x K(x)] \rightarrow [\exists x (K * \varphi)(x)]
\end{equation}
and (\ref{eq:agm2}) can be rewritten as
\begin{equation} \label{eq:agm4}
[\exists x ((K * \varphi)(x) \wedge \psi(x))]
\rightarrow [\forall y((K * (\varphi \wedge \psi))(y) \rightarrow ((K * \varphi)(y) \wedge \psi(y))].
\end{equation}

Postulates refer to a fixed knowledge base $K$, and are
implicitly universally quantified over formula symbols such as $\varphi, \psi$.
They express general requirements that are supposed to hold for all revising formulas.
Generalizing these examples, a postulate is defined as follows.

\begin{definition} (Postulate) \label{def:post}
A \emph{postulate} $P$ is a first-order sentence with unary predicate symbols $K, \varphi_1, \ldots, \varphi_\ell$ and $K * \mu_1, \ldots, K * \mu_m$,
where $\mu_1, \ldots, \mu_m$ are Boolean combinations of $\varphi_1, \ldots, \varphi_\ell$.

A revision operator satisfies a postulate for a knowledge base $K$ if the postulate holds for all $\varphi_1, \ldots, \varphi_\ell$,
with the variables ranging over the set of truth assignments.
\end{definition}

This definition covers all postulates in \cite{KM} and in Section 7.3 of \cite{H1999}.

\section{Characterizability}

As we consider partial preorders that are faithful for a knowledge base, we introduce the following property of partial preorders.

\begin{definition} (Regular partial preorders) \label{def:reg}
A partial preorder is \emph{regular} if
\begin{enumerate}
\item every minimal element is smaller than any non-minimal element

\item the number of elements is a power of 2.
\end{enumerate}
\end{definition}

An example of a non-regular partial preorder is the 4-element partial order with $a < b, c < d$ and no other comparability.
Condition 1 is satisfied, for example, if there is a unique minimal element.

\begin{definition} ({\cal R}-revision operator)
Let ${\cal R}$ be a family of regular partial preorders.
Let $K$ be a knowledge base and $*$ be a revision operator for $K$.
Then $*$ is an ${\cal R}$-revision operator iff there is a
faithful partial preorder $F = (R, t)$ for $K$, with $R \in {\cal R}$, representing $*$ using minimization.
\end{definition}

Using this definition, a formal definition is given of characterizability.

\begin{definition} (Characterization, characterizability)
A finite set of postulates ${\cal P}$ \emph{characterizes} ${\cal R}$-revision operators if
for every knowledge base $K$ and every revision operator $*$ for $K$ the following holds:
$*$ satisfies the postulates in ${\cal P}$ iff it is an ${\cal R}$-revision operator.

The family of ${\cal R}$-revision operators is \emph{characterizable} if there is a finite set of postulates characterizing ${\cal R}$-revision operators.
\end{definition}

It may be assumed \emph{w.l.o.g.} that ${\cal P}$ consists of a single postulate.


\section{A non-characterizable class}

A partial order is a \emph{crown} if it has elements $a_1, \ldots, a_s$ and $b_1, \ldots, b_s$ for some $s$, and the comparabilities are
$a_i > b_i$ and $a_i > b_{i+1}$ for every $i$, where the indices are meant cyclically.
A partial order is a \emph{double crown} if it consists of two crowns with pairwise incomparable elements.
An \emph{extended double crown} has additional elements that are smaller than all other elements.
An \emph{extended crown} is a crown with
additional elements that are smaller than all other elements.
Thus extended double crowns and extended crowns satisfy the first condition of Definition~\ref{def:reg}
and if their size is a power of 2 then they are regular.


Let ${\cal R}_0$ be the family of extended double crowns with size a power of 2,
and ${\cal R}_1$ be the family of extended crowns with size a power of 2.
The following theorem gives a sufficient condition for non-characterizability.
Its proof is given in the next two sections.


\begin{theorem} \label{th:main}
Let ${\cal R}$ be any family of regular partial orders containing ${\cal R}_0$ and disjoint from ${\cal R}_1$.
Then the family of ${\cal R}$-revision operators is not characterizable.
\end{theorem}

As a special case of the theorem, we formulate one specific non-characterizable family of revision operators.
A partial preorder $R$ is \emph{regular-disconnected} if it is regular, and the partial preorder obtained from $R$ by
removing its minimal elements has a disconnected comparability graph.

\begin{corollary}
Let ${\cal R}$ be the family of regular-disconnected partial orders.
Then the family of ${\cal R}$-revision operators is not characterizable.
\end{corollary}

The corollary follows directly from Theorem~\ref{th:main} and the definitions.

\section{Translation lemma}

Now we define a translation of postulates into sentences over an extension of the language of partial preorders.
The language of partial preorders contains a binary relation symbol $\le$ and equality.

The translated sentences also contain additional unary predicate symbols $A_1, \ldots, A_\ell$.
These correspond to propositional formulas $\varphi_1, \ldots, \varphi_\ell$ occurring in the postulates.
Given a Boolean combination $\mu$ of $\varphi_1, \ldots, \varphi_\ell$, we denote by $\hat{\mu}$ the first-order formula
obtained by replacing the $\varphi$'s with $A$'s. 
For instance, for
$\mu(x) = (\varphi_1 \wedge \varphi_2)(x)$ becomes $\hat{\mu}(x) = A_1(x) \wedge A_2(x)$.

Given a formula $\nu$ over the language $\le, A_1, \ldots, A_\ell$ with a single free variable $x$ we write $\min\nolimits_{\le}^{\nu}$ for
a formula expressing that $x$ is a minimal element satisfying $\nu$, i.e.,
\[ \min\nolimits_{\le}^{\nu}(x) \,\, \equiv \,\, \nu (x) \wedge \forall y (\nu(y) \rightarrow \neg (y < x)).\]
Minimal elements in the partial preorder are defined by
\[ \min\nolimits_{\le}(x) \,\, \equiv \,\, \forall y (\neg (y < x)).\]


\medskip

\begin{definition} (Translation) \label{def:tr}
The \emph{translation} $\tau(P)$ of a postulate $P$ is the sentence obtained from $P$ by replacing
\begin{enumerate}
\item every occurrence of $K(x)$ with $\min\nolimits_{\le}(x)$

\item every occurrence of $\varphi_i(x)$ and $\mu_i(x)$ with their ``hat" versions

\item every occurrence of $K * \mu_i$ with  $\min\nolimits_{\le}^{\hat{{\mu}_i}}(x)$
\end{enumerate}
\end{definition}

The translation is a first-order sentence over the predicate symbols $\le, A_1, \ldots, A_\ell$.

\medskip

\begin{example} (Translation of  postulate (\ref{eq:agm4})) Let us replace $\varphi$ and $\psi$ with $\varphi_1$ and $\varphi_2$ to be consistent
with the general notation:
\begin{equation} \nonumber
[\exists x ((K * \varphi_1)(x) \wedge \varphi_2(x))]
\rightarrow [\forall y((K * (\varphi_1 \wedge \varphi_2))(y) \rightarrow ((K * \varphi_1)(y) \wedge \varphi_2(y))].
\end{equation}
Applying Definition~\ref{def:tr} we get
\[ [\exists x (\min\nolimits_{\le}^{A_1}(x) \wedge A_2(x))]
\rightarrow [\forall y(\min\nolimits_{\le}^{A_1 \wedge A_2}(y) \rightarrow (\min\nolimits_{\le}^{A_1}(y) \wedge A_2(y))].\]

\end{example}

Given $K, \varphi_1, \ldots, \varphi_k$ and a faithful partial preorder $F$ for $K$, the $(\varphi_1, \ldots, \varphi_k)$-extension
of $F$ is determined in the standard way, by interpreting the unary predicate symbols $A_1, \ldots, A_k$ by $A_i(a) = \varphi_i(t(a))$.
The following proposition is a direct consequence of the definitions.

\begin{prop} \label{pr:tech}
Let $K$ be a knowledge base, $F = (R, t)$ be a faithful partial preorder for $K$ and let $*_F$ be the revision operator determined by $F$ using minimization.
Let $\varphi_1, \ldots, \varphi_\ell$ be propositional formulas and $P$ be a postulate.
Then $P$ is satisfied by $K$ for  $\varphi_1, \ldots, \varphi_\ell$ iff the $(\varphi_1, \ldots, \varphi_\ell)$-extension of $F$ satisfies $\tau(P)$.
\end{prop}

Now we formulate the connection to definability. In the formulation of we restrict ourselves to partial orders.



A \emph{universal monadic second-order sentence} is of the form
\[ \Phi = \forall A_1, \ldots, A_\ell \Psi, \]
 where $A_1, \ldots, A_\ell$ range over
unary predicates (or subsets) over the universe,
and $\Psi$ is a first-order sentence using the unary predicate symbols $A_1, \ldots, A_\ell$ in addition to the original language
(in our case $\le$ and equality). An existential second-order sentence is of the form $\Phi = \exists A_1, \ldots, A_\ell \Psi$.

\begin{definition} (Universal monadic second-order definability)
A family ${\cal R}$ of regular partial orders is universal monadic second-order definable if there is a universal monadic second-order sentence $\Phi$
such that for every regular partial order $R$ it holds that $R \in {\cal R}$ iff it satisfies $\Phi$.
\end{definition}

\begin{lemma} \label{lem:uni8} Let ${\cal R}$ be a family of regular partial preorders.
If the family of ${\cal R}$-revision operators is characterizable then ${\cal R}$ is universal monadic second-order definable.
\end{lemma}

\emph{Proof} Let ${\cal R}$ be a family of regular partial orders such that ${\cal R}$-revision operators are characterized by a postulate $P$.
We claim that ${\cal R}$ is defined by the universal monadic second-order sentence
\[  \Phi = \forall A_1, \ldots A_\ell \,\, \tau(P). \]

Assume that the regular partial order $R = (U, \le)$ is in ${\cal R}$. Let the number of its elements be $2^n$.
Let $t : U \to \{0, 1\}^n$ be an arbitrary bijection between $U$ and the set of truth assignments.
We get a faithful partial preorder $F = (R, t)$
for some knowledge base $K_n$, and thus the corresponding revision operator $*_F$ is an ${\cal R}$-revision operator.
Therefore $*_F$ satisfies $P$.
Consider arbitrary unary relations $A_1, \ldots, A_\ell$ over the elements.
Applying Proposition~\ref{pr:tech} to the propositional formulas $\varphi_1, \ldots, \varphi_\ell$ corresponding to $A_1, \ldots, A_\ell$,
it follows that $\le, A_1, \ldots, A_\ell$ satisfy $\tau(P)$.
Thus $R$ satisfies $\Phi$.

\medskip

Now assume that the regular partial order $R$ is not in ${\cal R}$.
Again, let $t : U \to \{0, 1\}^n$ be an arbitrary bijection between $U$ and the set of truth assignments.
We get a faithful partial preorder $F = (R, t)$ for some knowledge base $K_n$.
This determines a revision operator $*_F$.

We claim that $*_F$ is not an ${\cal R}$-revision operator. This follows if we show that $F$ is the only faithful partial order
determining $*_F$. Assume that for $F' = (R', t')$ with $R' = (U', \le')$ the revision operator $*_{F'}$ is the same.
Then, as revision operators are defined using minimization, for every pair of truth assignments $u, v$ it holds that
\begin{description}
\item $t^{-1}(u) < t^{-1}(v) \,\, \textrm{iff} \,\, K_n *_F \langle u, v \rangle = \langle u \rangle
\,\, \textrm{iff} \,\,  (t')^{-1}(u) < (t')^{-1}(v)$
\end{description}
and
\begin{description}
\item $t^{-1}(u) \sim t^{-1}(v)$  \,\, iff \,\, $K_n *_F \langle u, v \rangle = \langle u, v \rangle$
\,\, iff \,\,  $(t')^{-1}(u) \sim (t')^{-1}(v)$.
\end{description}
Thus $(t')^{-1} \circ t$ is an isomorphism from $R$ to $R'$.

Hence $*_F$ does not satisfy $P$.
So there are propositional formulas $\varphi_1, \ldots, \varphi_\ell$ such that the corresponding instance of $P$ is false.
By Proposition~\ref{pr:tech} the corresponding unary predicates $A_1, \ldots, A_\ell$
falsify $\tau(P)$. Hence $R$ falsifies $\Phi$. $\Box$

\section{Proof of Theorem~\ref{th:main}}

Two first-order structures are $q$-equivalent (denoted by $\equiv_q$) if they satisfy the same first-order sentences of quantifier rank at most $q$.
The $q$-round Ehrenfeucht - Fraiss\'e game over two structures is played by two players, Spoiler and Duplicator. In each round Spoiler picks
one of the structures and an element of that structure. Duplicator responds by picking an element in the other structure. After $q$ rounds Duplicator wins if the substructures
of picked elements in the two structures are isomorphic. Otherwise Spoiler wins. Duplicator has a winning strategy iff the two structures are $q$-equivalent.

The following result uses the notion of neighborhood in a general relational structure. We only use this result for undirected graphs with colored vertices,
where the $r$-neighborhood of a vertex $v$ is the set of vertices reachable from $v$ by paths of length at most $r$.

\begin{lemma} (Hanf-locality of first-order logic, see \cite{L}) \label{lem:hanf}
Let $\Psi$ be a first-order sentence with quantifier rank $q$, and let $r = (3^q - 1)/2$. Let $S_1$ and $S_2$ be two structures
with a bijection $f$ between their elements such that for every element $a$ of $S_1$ it holds that the $r$-neighborhoods of $a$ in $S_1$ and of $f(a)$
in $S_2$ are isomorphic. Then $S_1$ satisfies $\Psi$ iff $S_2$ satisfies it.
\end{lemma}

The proof of Theorem~\ref{th:main} is based on the undefinability of graph connectivity by existential monadic second-order sentences \cite{F, H}.
The basic fact is that every existential monadic second-order sentence satisfied by all cycles is also satisfied by some graph which is the union
of two cycles.
The following lemma is a variant of this result, based on the presentation in \cite{L}, with the required modifications.
The modifications are needed as a crown is somewhat different from a cycle and due to the presence of bottom elements Hanf-locality cannot be
applied directly.

\begin{lemma} \label{lem:mainle}
Let $\Phi$ be a universal monadic second order sentence over the language of $\le$ and equality. If every partial order in ${\cal R}_1$
falsifies $\Phi$ then some partial order in ${\cal R}_0$ also falsifies $\Psi$.
\end{lemma}

\emph{Proof}
Given an extended crown or extended double crown $M$, we define a structure $G_M$ over the language $E, L_1, L_2, L_3$, where $E$ is a binary relation and $L_1, L_2, L_3$
are unary relations. The ground sets of $M$ and $G_M$ are the same. Relation $E$ is the comparability relation of the crown or double crown, and
$L_1, L_2, L_3$ correspond to the maximal and minimal elements of the extended crown or extended double crown, resp., the bottom elements.
Thus $G_M$ is an undirected graph with vertices colored by the three colors $L_1, L_2, L_3$. We refer to $G_M$ as the colored graph of $M$.
The underlying undirected graph for extended crowns (resp., extended double crowns) is complete bipartite graph between a cycle
(resp. two cycles) and an independent set.
The structures $M$ and $G_M$ are inter-definable by simple first-order sentences.
In one direction
\[ (a < b) \equiv [E(a, b) \wedge ((L_2(a) \wedge L_1(b)) \vee (L_3(a) \wedge L_2(b)) \vee (L_3(a) \vee L_1(b)))].\]
 In the other direction
$E(a, b) \equiv (a < b) \vee (b < a)$, $L_3(a) \equiv min_\le (a)$ and $L_2, L_3$ can be defined similarly.

For the rest of the argument it is more convenient to switch to existential sentences.
We show that if an existential monadic second-order sentence
\[ \Phi = \exists A_1, \ldots, A_\ell \Psi \]
over the language $E, L_1, L_2, L_3$
is satisfied by the colored graph $G_{M_1}$ of every extended crown $M_1$ then it is also satisfied by the colored graph $G_{M_2}$ of some extended double crown $M_2$.
The lemma then follows directly.

Let $q$ be the quantifier rank of $\Psi$. Let $r = (3^q - 1)/2$ and $T = 2 \cdot 2^{\ell (2 r + 1)}$.

\begin{lemma} \label{lem:exis}
Let $M_1$ be an extended crown on at least $(4 r + 4) \cdot T$ elements.
For every $(A_1, \ldots, A_\ell)$-extension $G_{M_1}'$ of $G_{M_1}$ there is an
extended double crown $M_2$ of the same size and an $(A_1, \ldots, A_\ell)$-extension $G_{M_2}'$ of $G_{M_2}$
such that $G_{M_1}' \equiv_q G_{M_2}'$.
\end{lemma}

\emph{Proof of Lemma~\ref{lem:exis}} Let us consider the substructure $C_1$ of $G_{M_1}$ corresponding to vertices labeled $L_1, L_2$,
and its extension $(A_1, \ldots, A_\ell)$-extension $C_1'$. Thus $C_1$ is a 2-colored cycle.
The $r$-neighborhood of a vertex $a$ in $C_1$ is the set of vertices which can be reached from $a$ by a path of length at most $r$.
The neighborhood consists of two arcs of $r$ edges each.
The $r$-neighborhoods of vertices labeled $L_1$ (resp., $L_2$) in $C_1$ are isomorphic.
The extension adds an additional coloring with $2^\ell$ colors.
There are  $T$ isomorphism types of $r$-neighborhoods in $C_1'$.
Hence there are two elements $a$ and $b$ such that their distance on the cycle
is at least $2 r + 2$ and their $r$-neighborhoods are isomorphic.
Let $a'$ and $b'$ be the successors of $a$ and $b$ on the cycle (using some orientation).

Form $C_2'$ from $C_1'$ by deleting edges $(a, a')$ and $(b, b')$, and adding edges $(a, b')$ and $(b, a')$.
It follows from Lemma~\ref{lem:hanf} that $C_1' \equiv_q C_2'$.

Let $G_{M_2}'$ be the extended colored graph obtained from $C_2'$ by adding colored bottom vertices as in $G_{M_1}'$ .
Form $G_{M_2}$ from $G_{M_2}'$ by deleting the unary relations $A_1, \ldots, A_\ell$
and let $M_2$ be the extended double crown represented by $G_{M_2}$. The sizes of $M_1$ and $M_2$ are the same.

We claim that $G_{M_1}' \equiv_q G_{M_2}'$. It is sufficient to show that Duplicator wins the $q$-round first-order game on the
two structures. As $C_1' \equiv_q C_2'$, Duplicator wins the $q$-round first-order game on $C_1'$ and $C_2'$. Also, Duplicator
trivially wins the $q$-round first-order game on the two isomorphic sets of bottom vertices. As all edges are present between
the cycles and the bottom vertices, the combination of the two winning strategies gives a winning strategy on $G_{M_1}'$ and $G_{M_2}'$.
$\Box$ $\Box$

\bigskip

To conclude the proof of Theorem~\ref{th:main} assume that ${\cal R}$-revision operators are characterizable by a postulate $P$.
Then ${\cal R}$ is universal monadic second-order definable by Lemma~\ref{lem:uni8}.
Thus there is a universal mondaic second-order sentence satisfied by every partial order in ${\cal R}_1$ and falsified by every partial
order in ${\cal R}_0$, contradicting Lemma~\ref{lem:mainle}. $\Box$

\bibliographystyle{splncs}
 \bibliography{MSDb}

\end{document}